\documentclass[12pt,thmsa]{article}
\usepackage{graphicx}
\usepackage{amsmath}
\usepackage{amsthm}
\usepackage{a4wide}
\usepackage{mathrsfs}
\usepackage{relsize}
\usepackage{pdfpages}
\usepackage{amssymb}
\usepackage{hyperref}
\usepackage{graphicx}

\newcommand{\bR}{\mathbb{R}}

\theoremstyle{remark}

\title{Scalar Curvature, Flat Borromean Rings \\
and the 3-Body Problem}
\author{Michael Atiyah\\[1ex]
School of Mathematics\\
James Clerk Maxwell Building\\
The King's Buildings\\
Peter Guthrie Tait Road\\
Edinburgh EH9 3FD\\[4ex]
Dedicated to Isaac Newton and his Institute}

\date{\today}

\numberwithin{equation}{section}
\begin{document}
\maketitle

\section{History}

 The basic ideas of geometry go back to Euclid and Archimedes, followed later by Gauss, Minkowski and Weyl.  Topology also has ancient roots as evidenced by Alexander the Great cutting the Gordian Knot, and the Popes whose coat of arms became the Borromean rings.

 Isaac Newton, standing on the shoulders of Galileo and beyond Papal influence, used Greek ellipses to explain planetary motion and left to posterity the 3-body problem.

 The 3 Borromean rings of polished steel at the entrance to the Newton Institute in Cambridge\footnote{Image reproduced by kind permission of Peter Robinson, Bradshaw Foundation \copyright John Robinson/Edition Limit\'ee Paris 2007} provided me with inspiration, in August 2017, at the conference celebrating the 60th birthday of Simon Donaldson.
 Tantalized by Misha Gromov's lecture on scalar curvature, I decided to combine the two topics and relate them to Newton's planetary orbits.
 
$$\includegraphics[width=6cm]{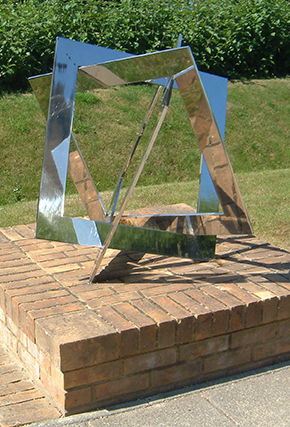}$$

 Visionary ideas are what we mathematicians grapple with, but illustrative special examples help the understanding.  They also lead us by analogy to more general and realistic examples.  The master whose style I try to emulate is the greatest mathematical physicist since Newton, James Clerk Maxwell.  He understood dynamical stability, with its topological underpinnings, more thoroughly than anyone else and I offer this short note as a tribute to him.

\section{The Topology and Geometry of Borromean Rings}

 The linking of circles in 3-space is familiar to us all and led to the formation of chains.  Gauss, Faraday and Maxwell explained electro-magnetism in terms of such links.

 The beauty and subtlety of the 3 Borromean rings is that each pair are unlinked but the triple holds together.  This mysterious property had, in medieval times, theological (even Trinitarian) overtones which is why it appealed to Popes, and would-be popes, like the Borromeo family.  Newton might have heard of this mystery but, as a secret unitarian in Trinity College, he would have kept off such dangerous ground.

 We all know that geometry and topology are siblings.  Geometry is about {\bf measurement} and topology is about {\bf shape}.  This family unity is embodied in the great theorem of Gauss, which starts locally with the area of spherical triangles, or polygons, and concludes globally with the formula for the Euler member of a closed surface as the average of the Gauss or scalar curvature (divided by $2\pi$).

 It is natural to ask if there is a similar formula relating the geometry and topology of Borromean rings.  For this to make sense the 3 rings must be geometrical, not just topological.  Each ring should be planar or spherical, with a boundary and the 3 rings should be placed in 3-space so that their mutual topology is Borromean.

 Looking at the sculptured steel outside the Newton Institute raises a challenging question, rather like the Rubik Cube or the model sailing boats in wine bottles.  How are they made or dismantled?  If you are the artist commissioned to make the sculpture, what precise instructions do you give the craftsman in the workshop?  For example, suppose you wanted to make the sculpture on the left of the Institute consisting of 3 squares of flat steel, each with a smaller square hole inside it?  The craftsman would need exact measurements and their tolerance : how accurate does he have to be?

 This sounds, and is, standard engineering design.  The sculptor, or his mathematical consultant, has to do "the sums" and instruct the craftsman, confident that it will give a stable work of art.

 Having done this once for the first steel sculpture at the Newton Institute, the design team then have to go through a similar procedure with the square replaced by a rhombus or a triangle, for the two sculptures on the right.

 As your sculptures attract attention the team might be flooded with orders, but each client would have different requirements.  Some might want polygons with many sides, perhaps infinitely many, so that polygons become ellipses (with Newtonian approval).   Clearly you will need an efficient algorithm to implement the process.  Does your firm have mathematical consultants of the necessary quality?  If they studied at Trinity or St. John's (the parent Colleges of the Newton Institute) then you might be in luck.

  As a Trinity man I feel it incumbent on me to try.  The next sections provide my mathematical report on these problems.

\section{The Iso-Perimetric Inequality}

 In the spirit of Maxwell I will consider a simple model based on the actual steel sculptures at Newton's Institute.  My Borromean rings will be modelled by flat regular polygons $C_i $ $(i=1,2,3)$ with $N$ sides ($N\geqslant 3$), having a border of width $\delta$.  See the figure for $N=4$ and one polygon $C$.
$$\includegraphics[width=10cm]{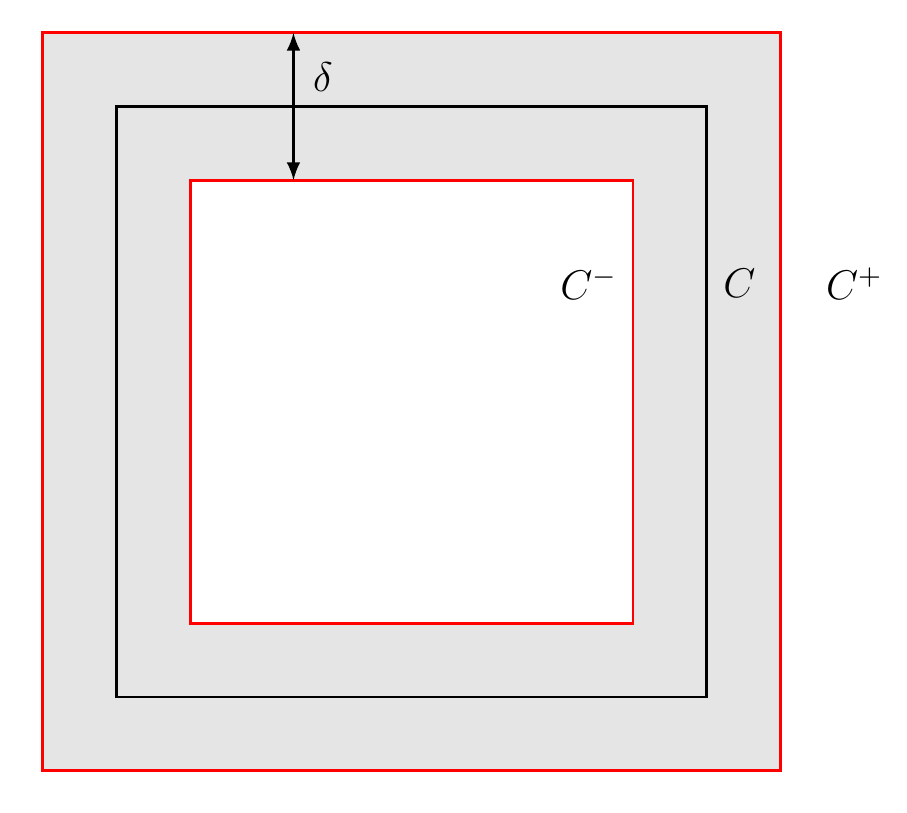}$$
The engineering problem is to take 3 such regular $N$-gons $C_1,C_2,C_3$ of the same width $\delta$ and then place them in 3-space (or on a 3-sphere) so that they form a Borromean triple.  Clearly the holes in the middle of the $N$-gons cannot be too small.  We look for an inequality which will give the critical size $d$, so that a Borromean triple of such (flat bordered) $N$-gons can be constructed provided $\delta<d$.

 Fortunately there is a famous inequality which is precisely what we need.  It is the {\bf iso-perimetric inequality}.  Moreover this inequality has been generalized over the centuries and goes by various names.  The simplest and most beautiful is due to Minkowski and concerns areas $A$ of plane {\bf convex} sets $P,Q$.  Minkowski considered the join of $P$ and $Q$, formed by all segments $\lambda p+(1-\lambda)q$ with $p\in P, q\in Q$ and $0\leqslant \lambda \leqslant 1$.  He discovered that the function
\begin{equation}
A(\lambda P+(1-\lambda) Q)
\end{equation}
 is a {\bf quadratic} in $\lambda$ whose first and last coefficients are $A(P)$ and $A(Q)$.  The middle term is called the Minkowski mixed area $A(P,Q)$.  The key result is that the quadratic (3.1) is {\bf indefinite}, i.e.
\begin{equation}
A(P,Q)^2 \geqslant 4 A(P) A(Q)
\end{equation}
 with equality only if $P=Q$  (up to translation and scale).

 The Isoperimetric Inequality follows by taking $Q$ to be a circle $S$ of (small) radius $\delta/2$  whose motion around the boundary $\partial Q$ of $Q$ sweeps  out a border of width $\delta$.  If $P$ is not convex, replacing it by its convex hull only sharpens the inequality.
 Note that the motion of the rolling ''marble" requires $\delta$ to be small relative to sharp bends of $\partial Q$.  But since the formula is exactly quadratic (not just approximately) we can start with a smooth boundary $\partial Q$, fix our $\delta$ and then let $\delta \rightarrow 0$.  Thus (3.2) holds for polygonal as well as smooth boundaries.

  There are deep relations between the Gauss Theorem on scalar curvature and Minkowski's Theorem on convex bodies.  They have both had extensive generalizations, in particular to higher dimensions.  They both sit at the crossroads of Geometry, Topology and Dynamics exemplified by the title of this paper.  But, as stated in the introduction, my aim in this short paper is to illustrate by well-chosen examples, not to develop general theory.

 The examples will start with specially symmetric ones, such as regular polygons, but since we are interested in topologically stable features we can deform to much more general configurations.

 Let us now return to our regular $N$-gon with vertices on a circle of radius $r$.  The formula for the area $A_N(r)$ is by elementary trigonometry

\begin{equation}
A_N(r)=\frac{r^2N}{2} \sin \frac{2\pi}{N}
\end{equation}

 Note that as $N\rightarrow \infty$ this converges to the formula $\pi r^2$ which is how Archimedes computed $\pi$.  As a function of $r$ it is just quadratic, as noted more generally by Minkowski, so that the area of the bordered $N$-gon is {\bf precisely} given by the first 2 non-zero terms of the Taylor expansion
\begin{equation}
A_N(r,\delta)=\frac{(2r\delta-\delta^2)N}{2} .\sin {\frac{2\pi}{N}}
\end{equation}
For both Gauss and Minkowski we can follow the formula by continuity as the polygon changes.  The difference is that, for Gauss, (3.4) would be an approximate formula only valid for small $\delta$ relative to the sharp bends on the boundary.  For Minkowski by contrast convexity makes (3.4) is exact and $\delta$ can be made uniformly small irrespective of the sharp bends in the boundary, when the quadratic ''error term" would lead to delta functions of boundary curvature (angles) for a polygon.

 Globally, for a surface of constant curvature, Minkowski and Gauss give the same (topological) answer.

 In this global form curvature can be both positive and negative, but the average has topological content.  This principle extends to all generalizations.  In the next section I will describe the pure topology of Borromean triples and show that it has hidden depths.

\section{The Topology of Borromean Rings}

  When topology became formalized as {\bf Algebraic Topology}, in the 20th Century, Borromean triples were seen to be related to other triples such as
\begin{itemize}
\item[(4.1)] the associator, that measures deviation from associativity, as in the octonions,
\item[(4.2)] the fundamental 3-form of a simple Lie group, starting with the unit quaternions,
\item[(4.3)] the vertices of a triangle,
\item[(4.4)] the Massey triple product in homotopy theory.
\end{itemize}
What is the right context of our example? There are many but the essential features are
\begin{itemize}
\item[(4.5)] 3 continuous closed curves in 3-space with the Borromean property, detected by (4.1) (4.2) or (4.4)
\item[(4.6)]  3 borders giving ribbon graphs
\item[(4.7)] 3 surfaces with various degrees of smoothness, for the interior and for the boundary
\item[(4.8)] 3 dimensions can represent static physical 3-space or dynamical 2-space incorporating time \{recall that Newton's planetary orbits are space-time graphs\}
\item[(4.9)] “Curve” in (4.5) can be interpreted as a compact set with “1-dimensional features” as in Gromov's talk.
\end{itemize}
These all relate to the way 2-planes sit in 3-space.  They can be usefully compared to (pairs) which lead to curvature, complex numbers and the way lines sit in 2-dimensions.

 The different aspects of Borromean triples listed above indicate that Borromean triples are not an amusing curiosity for popes or parties but are fundamental in 3-dimensional static phenomena (as is the M\"obius band).  This explains why stability is so important and how dynamics, which brings in time,  will relate to 4-dimensions as in Donaldson 60

 These general thoughts, when suitably developed, can be seen to underpin vast areas of natural philosophy and are much studied by scientists of all types, including mathematical ones.  But, since the entire universe in all its diversity is too vast to examine in full detail and at all scales, scientists long ago decided to study simple model examples.  This was fully understood by Maxwell.  As I explained I aim to be Maxwellian and pick simple illustrative examples.  That is why the Borromean steel sculptures at the Newton Institute are my inspiration.

 An example by itself is not enough.  We have to see it in a broader context and imagine ways, perhaps several, in which it can be generalized.  Then, and only then, does an example become an illustration.

 So, in this spirit, let me return to the geometry of the Newton Institute sculptures.

\section{A Geometric Result about Borromean Rings}

 The geometric data consists of 3 bordered regular $N$-gons with $N>3$ of width $\delta$ which are planar or spherical.  The figure in section 3 shows the planar case for $N=4$.  The question we ask is

\begin{itemize}
\item[(5.1)]  Can this data be assembled in 3-space ($\bR^3$ or $S^3$) by translation and scale changes to form a Borromean triple of $N$-gons?
\end{itemize}
The trigonometric calculations in Section 3 show that this is possible provided we have the inequality
\begin{itemize}
\item[(5.2)]    $A_N(r,\delta)  >  r\delta  \sin (2\pi/N)$
\end{itemize}
derived from the equality (3.4).  The "error" is the quadratic term 
\begin{itemize}
\item[(5.3)] $\delta^2/2.N \sin(2\pi /N)$
\end{itemize}

This answers the questions raised at the end of section 2. The desired bound $d$ depends on the scale set by the radius $r$ and, for regular $N$-gons $d=2r$ and we just get the inequality $\delta/r \leqslant 2$.

Embodied in planar steel, as at the Newton Institute, formula (5.2) gives the engineering team the precise measurements needed, for  any $N$, and the tolerance or error. As the width tends to zero the tolerance also tends to zero so that, in the (conformal) limit, steel is no longer a suitable material and the Borromean rings have to made of silk thread.

For irregular Borromean rings the exact formula (5.3) will be replaced by more
general ones, but the leading (topological) term will remain.

While (5.2) is elementary trigonometry, because of the symmetry, its natural context is that of Minkowski's mixed areas, or the Gauss scalar curvature. Generalizing (5.2) was therefore standard geometry over a century ago. In more recent times the topology and geometry of Borromean triples has become equally standard though more sophisticated.  Different experts develop different machinery and do not always use the same language.

For applied mathematicians, engineers or for Olympic pole vaulters, the position of the centre of gravity is crucial. Borromean geometry, based on the stability that comes from topology, may help win gold medals.  Maxwell would have understood it all and he would have explained how a few well-chosen examples are all you need. Since he is no longer here in person I have tried to stand in for him, like an actor who stands in for Charles Dickens.

    3 planets in elliptical orbits may have long term stability, held together by mysterious topology, if the parameters are right.  It seems to me more than likely that some versions of the famous 3-body problem will turn out to be Borromean,  perhaps in the rings of Saturn, studied by Maxwell or even in the Red Spot of Jupiter.  .

I hope I have whetted the appetite of the reader, and that the Newton Institute will continue to inspire natural philosophers for decades to come.

But let me conclude with two final sections extending the scope of the discussion.


\section{Force and Centres}

The Borromean rings at the Newton Institute were our starting point.  For artistic reasons they were highly symmetrical with an obvious centre.  But, as we start to deform them, their symmetry and centre get lost. How do we track the centre?  Newton found the answer by going back to Apollonius and conic sections.  The focus (one of a pair) was the appropriate centre and the eccentricity of the ellipse increased as the focus moved away from the centre. What Newton realized, and Apollonius did not, was that the focus was the centre of {\bf gravity}. In other words, gravitational force changes the geometry, and planetary orbits reflect this changed geometry. On the time scale of Earth-years the gravitational attraction of the sun does not change and planetary orbits are, with modest variations, stable for millennia.

Mathematicians, from Laplace onwards, wondered about stability in the “long-term”. Would the solar system survive forever? To an astronomer or an astrophysicist, this was an idealization too far.  A time-scale may not matter to the pure mathematician but, if it is not commensurate with physical time-scales, it is irrelevant: the “long-term” will never arrive.
Such metaphysical questions acquire more significance in two ways. On the one hand, if gravity were replaced by electro-magnetism, also based on the inverse-square law, planets could be replaced by electrons travelling in orbits at velocities close to the speed of light. This dramatic change of scales brings mathematical time-scales down to physical time-scales and long-term stability becomes important for engineers at CERN, as explained to me many years ago by J\"urgen Moser.

Another more fundamental change came, a century ago, with Einstein's Theory of General Relativity.  The inverse-square law is just the Newtonian linear approximation but non-linear GR has to incorporate Einsteinian corrections, even for the solar system, and the accuracy of GPS requires such corrections. If we also want to understand the magnetic fields that generate solar flares and affect our weather, then both kinds of corrections are necessary.

The LIGO triangles designed to detect gravitational waves were testing not just the statics of GR but also its dynamics, which brings in time, so Borromean rings might acquire a cosmic significance. This may not just be “$\pi$ in the SKY” since respected astronomers have been searching for “cosmic strings” both in theory and in observation.

The lessons I draw from this global (or cosmic) perspective are the following
\begin{itemize}
\item[(6.1)] the importance of forces in geometry
\item[(6.2)] the key role of the centres of force
\item[(6.3)] the role of centres in dynamical stability
\item[(6.4)] the importance of effective criteria for stability
\item[(6.5)] the understanding of different force scales
\item[(6.6)] the role of topology in qualitative understanding of statics and
         dynamics
\item[(6.7)] the role of geometry in quantitative versions of (6.6) with explicit tolerances.
\end{itemize}

Needless to say most of this is known to most scientists most of the time.  But I suspect that not all is known to all scientists all the time.

Perhaps it is time for all natural philosophers to get under the Newton Apple Tree and admire the Borromean Rings, just as our Druid ancestors admired the Art and Architecture of Stonehenge.

In the next and final section I will descend from Heaven to Earth.


\section{Down to Earth}

The 3 Borromean rings from which we started are the simplest geometric examples in 3-space.  They can be generalized in infinitely many ways and we are free to choose the model that fits our needs and our resources. There is no single multi-purpose model which will do everything.

I have already indicated that regular polygons can be replaced by convex polygons and that the parameters can all be allowed to vary.  For a bordered polygon with width $\delta$ we started with 3 polygons $C_i$ having the same width, but we can now allow different widths $\delta_i$. Ellipses appear when the polygon has an infinite number of sides, so we can end up with 3 ellipses of different widths sitting in 3-space.  Topologically they are either Borromean rings or non-Borromean rings, the Papal version of Hamlet's dilemma : {\it To $B$ or not to $B$}.

Mathematical techniques now exist to answer Hamlet. Shakespeare was an actor as well as a playwright, so he would have embraced the new technology.
The new technology discussed at Donaldson 60 would have solved Hamlet's dilemma, depriving the English language of its favourite quotation. All Hamlet had to do was to find out if the 3 rings could balance stably in the earth's gravitational field. The nearby apple tree would have shown him that the small convex set (the tolerance), defined by the widths $\delta_i$, should contain the centre of gravity. Fortunately the flat fenland of Cambridge would ensure that the gravitational field had essentially the same strength on both sides of the Institute.  Apples could therefore be substituted by steel rings, though the Borromean property would need a high speed camera to  capture any entanglement of  3 falling apples, the dynamics of which were famously observed by the occupant of Woolsthorpe manor. The engineers would not need to incorporate the relativistic Einsteinian corrections of section 6, since auroras are rarely seen at latitude $52^\circ$ North.

We thus end up with an effective numerical criterion for the 3 rings to be Borromean or non-Borromean. This describes all the steel sculptures that might be constructed.

As pointed out in section 3, Minkowski only dealt with convex bodies, which like spheres, have positive scalar curvature, while the Gauss formula allows regions of negative curvature. But passing to convex hulls swallows up the negative regions. There are many ways of understanding this difference.
\begin{itemize}
\item[(7.1)] A general Morse function $f$ on a closed surface will have maxima, minima and saddle points.
\item[(7.2)] The Hessian of $f$ at critical points distinguishes between the 3 cases in (7.1). It defines local geometries with scalar curvature positive at the maxima and minima, but negative at the saddle points.
\item[(7.3)] In the regions of positive curvature the Minkowski theory agrees with Gauss.
\item[(7.4)] The Gauss integral formula is about averages and Euler's topological invariant gives the dominant term.
\item[(7.5)] An average smooths out local irregularities, and gives good results when the averaging process is tailored to the nature of the irregularities.
\item[(7.6)] In all branches of science local irregularities occur and have to be averaged. Atoms in a molecule, planets in the solar system or galaxies in the cosmos are all “irregularities” in this sense.
\item[(7.7)] Analysts have developed very refined techniques of smoothing irregularities which lead to stochastic differential equations.
\item[(7.8)] Geometers have developed techniques of stability descending from Minkowski and used  in Mirror Symmetry.
\item[(7.9)] Physicists have techniques called renormalization   which provide powerful tools in quantum field theory.
\item[(7.10)] Number theorists since Euler and Riemann have used similar ideas to study the “irregular” distribution of primes, with the Riemann Hypothesis as the ultimate goal.
\end{itemize}

{\bf References}

There is a vast literature on geometry and dynamics, much too extensive to list here.  The great classic is of course Maxwell's Smith Prize Essay on Saturn's Rings. For more recent literature I can at least point to some  of my own papers which are relevant in various ways.

\begin{enumerate}
 \item M.F.Atiyah, {\it Convexity and commuting Hamiltonians}, Bull. Lond. Math. Soc. 14 (1981), 1--15.
\item --, {\it The Non-Existent Complex 6-Sphere}, arXiv:1610.09366v2 [math.DG] (2016)
\item --, {\it Geometric Models of Helium}, arXiv:1703.02532v1 [physics.gen-ph] (2017)
\item --, {\it Groups of Odd Order and Galois Theory}, to appear.
\item -- and J.Berndt, {\it Projective planes, Severi varieties and spheres}, Surveys in Differential Geometry VIII, Papers in Honor of Calabi, Lawson, Siu
and Uhlenbeck, International Press, Somerville, MA (2003).
\item --, R.Bott and L.G\aa rding, {\it Lacunas for hyperbolic differential operators with constant coefficients I.}, Acta Math. 124 (1970), 109--189.
\item -- and G.Segal, {\it Twisted K-theory and Cohomology}, in Inspired by S.S. Chern: A Memorial Volume in Honor of a Great Mathematician (Editor: P. A. Griffiths), Nankai Tracts in Mathematics Vol. 11, World Scientific
Publishing Co Inc. (2007).
\item J.C.Maxwell, {\it On the stability of the motion of Saturn's rings},
An essay which obtained the Adams Prize in 1856 in the University of Cambridge, Macmillan (1859)\\
Available online at \href{https://archive.org/details/onstabilityofmot00maxw}{https://archive.org/details/onstabilityofmot00maxw}
\end{enumerate}

\end{document}